\documentclass[10pt,a4paper]{article}
\usepackage[a4paper]{geometry}
\usepackage{amssymb,latexsym,amsmath,amsfonts,amsthm}
\usepackage{graphicx}

\usepackage{epsfig}

\newtheorem{Theorem}{Theorem}[section]
\newtheorem{Definition}[Theorem]{Definition}

\newtheorem{Corollary}[Theorem]{Collorary}
\newtheorem{Lemma}[Theorem]{Lemma}

\numberwithin{equation}{section}

\title{Higher order recurrences and row sequences of Hermite-Pad\'e approximation}

\date{\today}

\author{Guillermo L\'{o}pez Lagomasino, Yanely Zaldivar Gerpe\footnotemark[1]}%\footnotemark[1]}

\begin{document}

\maketitle
\renewcommand{\thefootnote}{\fnsymbol{footnote}}
\footnotetext[1]{Departamento de
Matem\'{a}ticas, Universidad Carlos III de Madrid, Avda. Universidad
30, 28911 Legan\'{e}s, Madrid, Spain. email: \{lago, yzaldiva\}\symbol{'100}math.uc3m.es.
Both authors received support from research grant MTM 2015-65888-C4-2-P of Ministerio de Econom\'{\i}a, Industria y Competitividad, Spain.}

\begin{abstract}
We obtain extensions of the Poincar\'e and Perron theorems for higher order recurrence relations and apply them to obtain an inverse type theorem for row sequences of (type II) Hermite-Pad\'e approximation of a vector of formal power series.

\end{abstract}

\medskip


\medskip


\maketitle

\section{Introduction}

\subsection{Background.}
Let $(f_n)_{n\geq 0}$ be a solution of the recurrence relation
\begin{equation}
\label{eq:1}
f_n + \alpha_{n,1}f_{n-1} + \cdots + \alpha_{n,m} f_{n-m} = 0, \qquad n \geq m,\end{equation}
with given initial conditions $f_0,\ldots,f_{m-1}$. If for all $n$ the coefficients $\alpha_{n,1},\ldots,\alpha_{n,m}$ are given, it is well known that the solution space of the recurrence relation is a vector space of dimension $m$.

Such recurrence relations appear and play a central role in many fields of mathematics: number theory, difference equations, continued fractions, and approximation theory to name a few. In the general theory, two results due to H. Poincar\'e \cite{Poin} and O. Perron \cite{Per1,Per2} single out (see also \cite{Gel}).

Assume that
\begin{equation}
\label{limit}
\lim_{n\to \infty} \alpha_{n,j} = \alpha_j, \qquad j=1,\ldots,m, \qquad \alpha_m \neq 0.
\end{equation}
Define the so called characteristic polynomial of \eqref{eq:1}
\begin{equation}
\label{char}
p(z) = z^m + \alpha_1 z^{m-1} + \cdots + \alpha_m = \prod_{j=1}^m (z - \lambda_j)
\end{equation}

\medskip

\noindent
{\bf Poincare's Theorem.} Suppose that $0 < |\lambda_1| < \cdots < |\lambda_m|$. Then, any solution $(f_n)_{n\geq 0}$ of \eqref{eq:1} verifies that either $f_n = 0, n\geq n_0,$ or
\begin{equation}
\label{ratio}
\lim_{n\to \infty} \frac{f_{n+1}}{f_n} = \lambda_k,
\end{equation}
where $\lambda_k$ is one of the roots of the characteristic polynomial.

\medskip

\noindent
{\bf Perron's Theorem.} Suppose that $0 < |\lambda_1| < \cdots < |\lambda_m|$ and $\alpha_{n,m} \neq 0, n \geq m$. Then, there exists a fundamental system of solutions $(f^{(k)}_{n})_{n\geq 0}, k=1,\ldots,m$ of \eqref{eq:1} such that
\begin{equation}
\label{ratio2}
\lim_{n\to \infty} \frac{f^{(k)}_{n+1}}{f^{(k)}_{n}} = \lambda_k, \qquad k=1,\ldots,m.
\end{equation}

\medskip

Each solution $(f_n)_{n\geq 0}$ of \eqref{eq:1} can be associated with a Taylor series. Namely, 
\begin{equation}
\label{Taylor}
f(z) = \sum_{n\geq 0} f_n z^n.
\end{equation}
In the sequel we will frequently identify a solution of \eqref{eq:1} and its associated Taylor series. The analytic properties of an analytic element are encoded in its Taylor coefficients. For example, under the assumptions of Perron's Theorem, if $f^{(k)}(z)$ denotes the Taylor series associated with $(f^{(k)}_{n})_{n\geq 0}$ from \eqref{ratio2} it immediately follows that the radius of convergence  $R_0(f^{(k)})$ of $f^{(k)}$ equals $1/|\lambda_k|$. (Do not confuse $f^{(k)}$ with a derivative.) Thanks to a very deep result of E. Fabry \cite{fabry}, more can be said.

\medskip

\noindent {\bf Fabry's Theorem.} Given a Taylor series $f$ whose coefficients verify \eqref{ratio} we have that $R_0(f) = |\lambda_k|^{-1}$ and $\lambda_k^{-1}$ is a singular point of $f$.

\medskip

The object of this paper is to extend the results of Poincar\'e and Perron to more general recurrence relations for which the zeros of the characteristic polynomial do not necessarily have different absolute value, and describe some of the analytic properties of the functions associated  with a system of fundamental solutions of the recurrence relation.

\medskip

Set
\begin{equation}
\label{alphan}
\alpha_n(z) := 1 + \alpha_{n,1}z + \cdots + \alpha_{n,m}z^m.
\end{equation}
Unless otherwise stated, in the sequel we will assume that $\alpha_{n,m} \neq 0, n \geq m$.
If $[f]_n$ denotes the $n$-th Taylor coefficient of a formal power series $f$, then \eqref{eq:1} adopts the form
\begin{equation}
\label{eq:2}
[f\alpha_n]_n = 0, \qquad n \geq m.
\end{equation}
Let
\[ {\mathcal{P}}_{n} = \{\zeta_{n,1},\ldots,\zeta_{n,m}\},
\quad n\ge m,
\]
denote the  collection of zeros of $\alpha_{n}$ repeated according to
their multiplicity. Set
$$
S=\sup_{N\ge m}\inf_{n \geq N}\left\{|\zeta_{n,k}|:  \zeta_{n,k} \in  {\mathcal{P}}_{n} \right\}
$$
and
$$
G=\inf_{N\ge m}\sup_{n \geq N}\left\{|\zeta_{n,k}|:   \zeta_{n,k} \in  {\mathcal{P}}_{n}\right\}.
$$

The following result is a consequence of \cite[Theorem 2.2]{cacoq2}

\begin{Theorem} Assume that $S > 0$ and $G < +\infty$. Then, any nontrivial solution $f$ of \eqref{eq:1} verifies $0 < c \leq R_0(f) \leq C < +\infty$, where $c$ and $C$ only depend on $S$ and $G$.
\label{teo:1}
\end{Theorem}

This means that under such general conditions, every solution of \eqref{eq:1} has a singularity in the annulus $\{z: c \leq |z| \leq C\}$. A more precise result may be given when the sequence of polynomials $(\alpha_n)_{n \geq m},$ has a limit.

\medskip

In the following, we assume that \eqref{limit} takes place and
\[\lim_{n \to \infty} \alpha_n(z) := \alpha(z) = \prod_{k=1}^m \left(1 - \frac{z}{\zeta_k}\right) ,\qquad \deg (\alpha) = m.\]
According to \eqref{char}, we have $\alpha(z) = z^{m}p(1/z)$; therefore, the zeros of $\alpha$ and of the characteristic polynomial $p$ are reciprocals one of the other.

\medskip

Putting together Theorems 1 and 2 of V. I Buslaev in  \cite{Bus1}, we can formulate the following:

\medskip

\noindent
{\bf Buslaev's Theorem.} Assume that \eqref{limit} takes place and $f$ is a non-trivial solution of \eqref{eq:1}. Then $R_0(f)$ is equal to the absolute value of one of the zeros of $\alpha$, the coefficients of $f$ satisfy a reduced recurrence relation of the form
\begin{equation}
\label{cociente}
f_n + \beta_{n,1}f_{n-1} + \cdots + \beta_{n,\ell} f_{n-\ell} = 0,\qquad \lim_n \beta_{n,k} = \beta_k, \qquad k=1,\ldots,\ell,
\end{equation}
where $\ell (\leq m)$ is equal to the number of zeros of $\alpha$ on the circle $\{z: |z| = R_0(f)\}$, all the zeros of $\beta(z) = 1+\beta_1z + \cdots\beta_\ell z^{\ell}$ lie on that circle, $\beta$ divides $\alpha$, and at least one of its zeros is a singular point of $f$.

\medskip

When all the zeros of $\alpha$ have distinct absolute value, it is easy to see that Buslaev's theorem reduces to Poincare's result. A natural question arises, is there a fundamental system of solutions of \eqref{eq:1} such that every zero of $\alpha$ is a singular point of at least one function in the system?

\medskip

The answer to this question is positive if it is known in advance that \eqref{eq:1} has a fundamental system of solutions of the form $\{f,zf,\ldots, z^{m-1}f\}$ for some (formal) Taylor expansion $f$ about the origin. This is a consequence of a result of S.P. Suetin \cite[Theorem 1]{suetin1985} where he solves a conjecture posed by A. A. Gonchar in \cite{gon2} in the context of inverse type results for row sequences of Pad\'e approximants. In the last section, we will see the connection between the study of row sequences of Pad\'e (and Hermite-Pad\'e) approximation and the study of the general theory of recurrence relations such as \eqref{eq:1}. In \cite[Theorem 1]{gon2}, again within the study of row sequences of Pad\'e approximation,  Gonchar studied the case when
\begin{equation}
\label{rate}
\limsup_{n\to \infty} \|\alpha_n - \alpha\|^{1/n} = \theta < 1;
\end{equation}
that is, when the zeros of the polynomials $\alpha_n$ converge with geometric rate to the zeros of $\alpha$. In \eqref{rate}, $\|\cdot\|$ represents (for example) the coefficient norm in the finite dimensional space of polynomials of degree $\leq m$. In this situation, it may be proved  that $f$ has in the disk of radius $\max\{|\zeta_k|: k=1,\ldots, m\}/\theta$ exactly $m$ poles which coincide, taking account of their order, with the zeros of $\alpha$ taking  account of their multiplicity. Additionally, $\max\{|\zeta_k|: k=1,\ldots, m\}/\theta$ is the radius of the largest disk centered at the origin to which $f$ has a meromorphic extension with no more than $m$ poles. Gonchar's theorem was extended to the context of row sequences of (type II) Hermite-Pad\'e approximation in \cite[Theorem 1.4]{cacoq2}. As we shall see in the final section, the study of row sequences of Hermite-Pad\'e approximation is equivalent to the study of general solutions of recurrence relations of type \eqref{eq:1}. Having this in mind, we will present an extension of the Poincar\'e and Perron theorems to general recurrence relations \eqref{eq:1} verifying \eqref{limit} and some zeros of the polynomials  $\alpha_n$ do not converge geometrically.

\medskip

\subsection{Statement of the main result.} Without loss of generality, let the zeros of $\alpha$ be enumerated in such a way that
\begin{equation}
\label{enumeration}
 0 < |\zeta_1| \leq |\zeta_2| \leq \cdots \leq |\zeta_m|.\end{equation}
If several $\zeta_k$ coincide, they are enumerated consecutively. We shall also assume that the zeros in the collections of points $\mathcal{P}_n$ are indexed so that
\begin{equation}
\label{limit2}
\lim_{n\to \infty} \zeta_{n,k} = \zeta_k, \qquad k=1,\ldots,m.
\end{equation}
Several circles centered at the origin may contain more than one zero of $\alpha$ (or a zero of $\alpha$ of multiplicity greater than $1$). Let $C$ be one  such circle (if any). Let $\zeta_j,\zeta_{j+1},\ldots,\zeta_{j+N-1}$ be the zeros of $\alpha$ lying on $C$. In this case, we assume that
\begin{equation}
\label{several}
\limsup_{n\to \infty} |\zeta_{n,k} - \zeta_k|^{1/n} < 1, \qquad k=j,\ldots,j+N-1.
\end{equation}
The existence of such a circle $C$ is not required. If such a circle does not exist, all the zeros of $\alpha$ have distinct absolute value and we are in the situation of the Poincar\'e and Perron theorems.

\begin{Theorem}
\label{teo:2}
Assume that \eqref{limit2} and \eqref{several} take place. Then there is a fundamental system of solutions $\{f^{(1)},\ldots,f^{(m)}\}$ of \eqref{eq:1} such that  $R_0(f^{(k)}) = |\zeta_k|, k=1,\ldots,m,$  and $\zeta_k$ is a singular point of $f^{(k)}$. Each $\zeta_k$ verifying \eqref{several} is a pole of $f^{(k)}$. Moreover, if $\zeta_k$ is a zero of multiplicity $\tau$ and $\zeta_k = \zeta_{k+1} = \cdots = \zeta_{k+\tau-1}$ then for each $s= 1,\ldots,\tau$, $f^{(k + s-1)}$ is analytic in a disk of radius larger than $|\zeta_k|$ except for a pole of exact order $s$ at $\zeta_k$.
\end{Theorem}

\medskip

If the zeros of $\alpha$ have distinct absolute value, the statement of Theorem \ref{teo:2} is deduced directly from the Perron and Fabry theorems. When \eqref{several} takes place for $k=1,\ldots,m$, the thesis of Theorem \ref{teo:2} is a consequence of $(b) \Rightarrow (a)$ in \cite[Theorem 1.4]{cacoq2}.

\section{Proof of Theorems \ref{teo:1} and \ref{teo:2}.}

\subsection{Proof of Theorem \ref{teo:1}.} This result is a particular case of \cite[Theorem 2.2]{cacoq2}. To see this, for the benefit of the reader, we establish a connection between the notation employed here and in that paper.

Let $f$ be any non-trivial solution of \eqref{eq:1} and $t_n = t_n(f)$ the polynomial part of degree $\leq n-1$ of $\alpha_n f$. Due to \eqref{eq:1} and \eqref{eq:2} it follows that
\[(\alpha_n f - t_n)(z) = \mathcal{O}\left({z^{n+1}}\right).\]
This means that the rational function $t_n/\alpha_n$ is what in \cite{cacoq2} is called an $(n,m,m^*)$ incomplete Pad\'e approximation of $f$ with $m^* = 1$.

In \cite{cacoq2}, $\alpha_n$ is denoted $q_{n,m}$, $t_n$ is denoted $p_{n,m}$, $\lambda_n$ is the degree of the common zero which $p_{n,m}$ and $q_{n,m}$ may have at $z=0$, $m_n = \deg(q_{n,m})$, and $\tau_n = \min\{n-m^*-\lambda_n - \deg(p_{n,m}),m-\lambda_n - m_n\}$. In our case, since $\alpha_n(0) = 1$ and $\deg(q_{n,m}) = m$, we have that $\lambda_n = 0, m_n = m, \tau_n = 0,$ and $m^* = 1$. Thus, the assumptions in (i) and (ii) of \cite[Theorem 1.2]{cacoq2} are fulfilled and, therefore, $0 < R_0(f) < \infty$ as claimed. The proof gives lower and upper estimates of $R_0(f)$ depending on $S$ and $G$  which imply the last assertion of Theorem \ref{teo:1}.               \hfill $\Box$

\subsection{Comments on Buslaev's Theorem.} In \cite[Theorem 1]{Bus1} it is required that the solution considered is not a polynomial. We have stipulated that $\alpha_{n,m} \neq 0, n \geq m$. This restriction implies that any non-trivial solution cannot be a polynomial. In fact, the contrary would imply that $f_n = 0$ for all $n \geq n_0$ and solving the recurrence backwards (taking advantage that $\alpha_{n,m} \neq 0, n \geq m,$) we conclude that $f_n = 0, n\geq 0,$ and the solution would be the trivial one.
\medskip

In \cite[Theorem 2]{Bus1} the author imposes that $0 < R_0(f) < \infty$. Our assumption imply this as Theorem \ref{teo:1} shows. The polynomials that we have denoted $\alpha_n$ play the role of the functions $\alpha_n$ in \cite[Theorem 2]{Bus1}.

\subsection{Some auxiliary results.}

The proof of Theorem \ref{teo:2} is somewhat constructive. We start out from a fundamental system of solutions of \eqref{eq:1} and through analytic continuation, carried out in successive steps, we find another fundamental system of solutions which fulfills the desired properties. As we carry out these steps, we find collections of solutions which according to Buslaev's theorem have radius of convergence equal to the absolute value of a zero of $\alpha$. On any such  circle, there may fall one or several zeros of $\alpha$. The proof distinguishes two cases. The first when all the zeros on the circle satisfy \eqref{several}. In this case, the analytic continuation is based on \cite[Corollary 2]{Bus1}. If the circle contains only one zero of $\alpha$ of multiplicity 1 and we do not have \eqref{several}, we adapt a proof of Perron's Theorem given by M.A. Evgrafov in \cite{Evg} to continue the process.

\medskip

Here, we state \cite[Corollary 2]{Bus1} in the form of a lemma for the reader's convenience. We wish to mention that \cite[Theorem 2.6]{cacoq2} plays the same role when all the zeros satisfy \eqref{several}.

\begin{Lemma}
\label{lem:1}
Suppose that the assumptions of Buslaev's theorem hold and $f$ is a non trivial solution of \eqref{eq:1}. Let $\zeta_j,\ldots,\zeta_{j+N-1}, N \geq 1,$ be the zeros of $\alpha$ on the circle $\{z: |z| = R_0(f)\}$ and \eqref{several} takes place. Then $R_0(g) > R_0(f)$, where $g(z) = \prod_{k = j}^{j+N-1}(z - \zeta_k) f(z)$
\end{Lemma}

Evgrafov's proof of Perron's theorem is based on several lemmas. His paper has not been translated so for completeness we include the proof of those lemmas that we will use. The next two lemmas are identical to his statements. Lemma \ref{lem:4} has been adapted to cover our more general situation.

\begin{Lemma}
\label{lem:2}
Let $(f_n)_{n\geq 0}$ be a solution of \eqref{eq:1} and let $(\gamma_n)_{n \geq 0}$ be a sequence such that $\gamma_n \neq 0, n\geq 0$. Then $(F_n = f_n/\gamma_n)_{n\geq 0}$ is a solution of the recurrence relation
\begin{equation}
\label{eq:3}
 F_n + \alpha_{n,1}'F_{n-1} + \cdots + \alpha_{n,m}'F_{n-m} = 0, \qquad n \geq m,\end{equation}
where
\[\alpha_{n,j}' = \alpha_{n,j} \frac{\gamma_{n-j}}{\gamma_{n}}, \qquad j=1,\ldots,m.\]
Moreover, if $\,\lim_{n\to \infty} \gamma_{n+1}/\gamma_n = 1$ and the $\alpha_{n,j}, j=1,\ldots,m,$ verify \eqref{limit} then so do the $\alpha_{n,j}', j=1,\ldots,m $ and the recurrences \eqref{eq:1} and \eqref{eq:3} have the same characteristic polynomial.
\end{Lemma}

\noindent {\bf Proof.} Indeed, if we divide \eqref{eq:1} by $\gamma_n$, we obtain that for all $n\geq m$
\[0= \frac{f_n}{\gamma_n} + \alpha_{n,1}\frac{\gamma_{n-1}}{\gamma_n}\frac{f_{n-1}}{\gamma_{n-1}} + \cdots + \alpha_{n,m} \frac{\gamma_{n-m}}{\gamma_n}\frac{f_{n-m}}{\gamma_{n-m}},\]
which is \eqref{eq:3}.

Now, $\lim_{n\to \infty} \gamma_{n+1}/\gamma_n = 1$ implies that $\lim_{n\to \infty} \gamma_{n-j}/\gamma_n = 1, j=1,\ldots,m$, so  from \eqref{limit}  it follows that
\[\lim_{n\to \infty} \alpha_{n,j}' = \lim_{n\to \infty} \alpha_{n,j} = \alpha_j,\qquad j=1,\ldots,m.\]
Therefore, in that case both equations have the same characteristic polynomial. \hfill $\Box$

\medskip

It is easy to check that $(f_n^{(j)})_{n\geq 0}$, $j=1,\dots,N, 1\leq N \leq m,$ are linearly independent solutions of \eqref{eq:1} if and only if
$(f_n^{(j)}/\gamma_n)_{n\geq 0}$, $j=1,\dots, N,$ (where $\gamma_n\neq 0$, $n\geq 0$) are linearly independent solutions of \eqref{eq:3}.

\begin{Lemma}
\label{lem:3}
If equation \eqref{eq:1} has the solution $(\lambda^n)_{n\geq 0}, \lambda \neq 0$, then it's left hand side may be expressed in the form
\begin{equation}
\label{eq:4}
F_{n-1} + \beta_{n,1} F_{n-2} + \cdots + \beta_{n,m-1}F_{n-m} = 0
\end{equation}
 where
 \[F_n = f_{n+1} - \lambda f_n.\]
Additionally, if the $\alpha_{n,j}$ verify \eqref{limit} then
\begin{equation}
\label{limit*}
\lim_{n,j} \beta_{n,j} = \beta_j, \qquad j = 1,\ldots,m-1, \qquad \beta_{m-1} \neq 0
\end{equation}
and the characteristic polynomials of \eqref{eq:1} and \eqref{eq:4} are connected by the relation
\begin{equation}
\label{conexion}
(z - \lambda)(z^{m-1} + \beta_1 z^{m-2} + \cdots + \beta_{m-1}) = (z^{m} + \alpha_1 z^{m-1} + \cdots + \alpha_{m})
\end{equation}
\end{Lemma}

\noindent {\bf Proof.} Consider the polynomial
\[p_n(z) = z^m + \alpha_{n,1}z^{m-1}+ \cdots + \alpha_{n,m}.\]
Substituting the solution $(\lambda^n)_{n\geq 0}$ in \eqref{eq:1} and factoring out $\lambda^{n-m}$, we get
\[p_n(\lambda) = \lambda^m + \alpha_{n,1}\lambda^{m-1}+ \cdots + \alpha_{n,m} = 0;\]
therefore,
\[  \frac{p_n(z)}{z -\lambda} = z^{m-1} + \beta_{n,1} z^{m-2} + \cdots + \beta_{n,m-1}. \]
That is
\begin{equation}
\label{conexion*} p_n(z) = (z-\lambda)(z^{m-1} + \beta_{n,1} z^{m-2} + \cdots + \beta_{n,m-1}).
\end{equation}
Equating coefficients of equal power of $z$, we obtain
\begin{equation}
\label{relacion}
\beta_{n,j}-\lambda \beta_{n,j-1} = \alpha_{n,j}, \qquad j=1,\ldots,m, \qquad \beta_{n,0} = 1, \qquad \beta_{n,m} = 0.
\end{equation}
In particular,
\[\beta_{n,m-1} = - \alpha_{n,m}/\lambda \neq 0.\]
From \eqref{relacion} the existence of $\lim_{n\to \infty} \alpha_{n,j}$ and $\lim_{n\to \infty} \beta_{n,j-1}$ imply the existence of
$\lim_{n\to \infty} \beta_{n,j}.$ Since $\beta_{n,0} = 1$, it follows that \eqref{limit} implies \eqref{limit*} and then \eqref{conexion} is immediate taking limit over $n$ in \eqref{conexion*}.

\medskip

It remains to verify \eqref{eq:4}. Put $F_n = f_{n+1} - \lambda f_n$ in the left hand side of \eqref{eq:4}. Using \eqref{relacion} and \eqref{eq:1}, we get
\[F_{n-1} + \beta_{n,1} F_{n-2} + \cdots + \beta_{n,m-1}F_{n-m} = \]
\[f_{n} - \lambda f_{n-1} + \beta_{n,1}(f_{n-1} - \lambda f_{n-2}) + \cdots + \beta_{n,m-1}(f_{n-m+1} - \lambda f_{n-m}) = \]
\[ f_n + (\beta_{n,1} - \lambda \beta_{n,0}) f_{n-1} + \cdots + (\beta_{n,m} -\lambda \beta_{n,m-1})f_{n-m} = 0\]
as we needed to prove. \hfill $\Box$

\medskip

Assuming that $(\lambda^n)_{n\geq 0}, \lambda \neq 0,$ is a solution of \eqref{eq:1}, it is easy to verify that $(f_n^{(j)})_{n\geq 0}$, $j=1,\ldots, N (\leq m-1)$ and $(\lambda^n)_{n\geq 0}$ constitute a system of linearly independent solutions of \eqref{eq:1} if and only if $(F_n^{(j)})_{n\geq 0}$, $j=1,\ldots,N,$ is a system of linearly independent solutions of \eqref{eq:4}, where $F_n^{(j)}=f_{n+1}^{(j)} -\lambda f_n^{(j)}, j=1,\ldots,N$.

\begin{Lemma}
\label{lem:4}
Suppose that $(F_n)_{n\geq 0}$ is such that $\limsup_{n\to \infty} |F_n|^{1/n} = \mu, \mu \neq |\lambda|$. Then, there exists a solution $(f_n)_{n\geq 0}$ of the equations $F_n = f_{n+1} - \lambda f_n$ such that
$\limsup_{n\to \infty} |f_n|^{1/n} = \mu$.
\end{Lemma}

\noindent {\bf Proof.} We give two different expression for the solution $(f_n)_{n\geq 0}$ depending on whether $\mu < |\lambda|$ or $\mu > |\lambda|$.
In the first case we set
\begin{equation}
\label{sol1}
f_n = - \frac{F_n}{\lambda} - \frac{F_{n+1}}{\lambda^2} - \frac{F_{n+2}}{\lambda^3} - \cdots ,
\end{equation}
and in the second
\begin{equation}
\label{sol2}
f_n = F_{n-1} + \lambda F_{n-2} + \cdots + \lambda^{n-1} F_0.
\end{equation}
We will see in a minute that \eqref{sol1} is convergent for each $n$. With this in mind, it is easy to verify that so defined the sequence $(f_n)_{n\geq 0}$ satisfies the required equations.

Let us verify that the numbers $f_n$ in \eqref{sol1} are finite and $\limsup_{n\to \infty} |f_n|^{1/n} = \mu$. Indeed, take $\varepsilon > 0$ such that $\mu + \varepsilon < |\lambda|$. From the assumption on the $F_n$ we get that for some constant $C \geq 1$
\[|F_n| \leq C(\mu + \varepsilon)^n, \qquad n \geq 0.\]
Consequently,
\begin{equation}
\label{acota}
|f_n| \leq C\left(\frac{(\mu + \varepsilon)^n}{|\lambda|} + \frac{(\mu + \varepsilon)^{n+1}}{|\lambda|^2} + \frac{(\mu + \varepsilon)^{n+2}}{|\lambda|^3} + \cdots \right)  = \frac{C (\mu + \varepsilon)^n }{|\lambda| - (\mu + \varepsilon)} < \infty.
\end{equation}
Moreover,  \eqref{acota} implies that
\[\limsup_{n\to \infty }|f_n|^{1/n} \leq \mu + \varepsilon.\]
Letting $\varepsilon$ tend to zero we get that $\limsup_{n\to \infty }|f_n|^{1/n} \leq \mu $.  On the other hand, $|F_n| \leq |f_{n+1}| + |\lambda f_n|$. This in turn implies that $\mu = \limsup_{n\to \infty} |F_n|^{1/n} \leq \limsup_{n\to \infty} |f_n|^{1/n}$, so indeed we have equality.

When $\mu > |\lambda|$, we proceed analogously. Using \eqref{sol2}, we have
\begin{equation}
\label{acota2}
|f_n| \leq C (\mu + \varepsilon)^{n-1} \left(1 + \frac{|\lambda|}{\mu + \varepsilon} + \cdots + \frac{|\lambda|^{n-1}}{(\mu + \varepsilon)^{n-1}}  \right) \leq \frac{C (\mu + \varepsilon)^{n}}{\mu + \varepsilon - |\lambda|}.
\end{equation}
From \eqref{acota2} we get that $\limsup_{n\to \infty} |f_n|^{1/n} \leq \mu$, and equality is derived just as before. \hfill $\Box$

\section{Proof of Theorem \ref{teo:2}.}

\noindent
Let ${\bf f}^1 = ( {f}^{1,1}, {f}^{1,2},  \ldots, {f}^{1,m})$  be a fundamental system of solutions of \eqref{eq:1} and $\zeta_1,\ldots,\zeta_m$ the collection of zeros of $\alpha$ enumerated as in \eqref{enumeration}. According to Buslaev's theorem
\[0 < |\zeta_1| \leq R_0( {f}^{1,j}) \leq |\zeta_m| < \infty, \qquad j=1,\ldots,m.\]
Let $D_0({\bf f}^1)$ denote the intersection of all the disks centered at the origin and radii  $R_0( {f}^{1,j}), j= 1,\ldots,m$ whose boundary we denote $\mathcal{C}_1$. By Buslaev's theorem, the circle $\mathcal{C}_1$ of radius $R_1:= R_0({\bf f}^1)$ contains at least one zero of $\alpha$. Moreover, several (and at least one) of the functions in ${\bf f}^1$ has radius of convergence equal to $R_1 = R_0({\bf f}^1)$. Let $\zeta_j,\ldots,\zeta_{j+N-1}$ be the collection of all the zeros of $\alpha$ lying on $\mathcal{C}_1$. We  distinguish the cases when $N\geq 2$ and when $N = 1$. We remark that so far we cannot assert that $|\zeta_j| = |\zeta_1|$; in principle, it may be larger.

\medskip

Suppose that $N \geq 2$. In this case, according to our assumptions, \eqref{several} takes place. Due to Lemma \ref{lem:1}, $R_0(\prod_{k=j}^{j+N}(z - \zeta_k) {f}^{1,\ell}) > R_0({\bf f}^1)= R_1, \ell=1,\ldots,m$; that is, either $ {f}^{1,\ell}$ has radius of convergence larger than $R_1$ to start with or, $ {f}^{1,\ell}$ has at most poles on  $\mathcal{C}_1$ at zeros of $\alpha$ and their order is less than or equal to the multiplicity of the corresponding zero of $\alpha$.

\medskip

Let us  find coefficients
$c_1,\ldots,c_{ m}$ such that
\begin{equation} \label{lasg} \sum_{\ell=1}^{ m} c_\ell  {f}^{1,\ell}
\end{equation}
is analytic in a neighborhood of $\overline{D_0(\mathbf{f}^1)}$.
Finding the coefficients $c_\ell$ reduces to solving a homogeneous linear
system of $N$ equations on $m$ unknowns. In
fact, if $\zeta$ is one of the zeros of $\alpha$ on  $\mathcal{C}_1$ and it has multiplicity $\tau$ we obtain $\tau$ equations choosing the
coefficients $c_\ell$ so that
\begin{equation} \label{ecuacion} \int_{|\omega - \zeta| =
\delta} (\omega - \zeta)^\nu \left(\sum_{\ell=1}^{\bf m} c_\ell
 {f}^{1,\ell}(\omega)\right) d\omega = 0, \qquad \nu=0,\ldots,\tau -1.
\end{equation}
where $\delta$ is sufficiently small.  We do the same with each
distinct zero of $\alpha$ on  $\mathcal{C}_1$. The homogeneous linear system of $N$ equations so
obtained has $m - N_1, N_1 \leq N,$ linearly independent solutions, where $N_1$ equals the rank of the linear system of equations.
Denote the solutions of the linear system by ${\bf c}^{1,j},$ $ j=N_1 + 1,\ldots,m $.

\medskip

Set
\[ {\bf c}^{1,j} = (c_{1}^{1,j},\ldots,c_{m}^{1,j}),
\qquad j=N_1 +1,\dots,m ,
\]
and
\[ {f}^{2,j} = \sum_{\nu=1}^m c_{\nu}^{1,j}  {f}^{1,\nu}, \qquad j=N_1+1,\ldots,m.\]
We wish to emphasize several points:
\begin{enumerate}
\label{properties}
\item The collection of functions ${\bf f}^2 = ( {f}^{2,N_1 +1},\ldots, {f}^{2,{m }})$ is made up of nontrivial linearly independent solutions of \eqref{eq:1}.
\item Because of \eqref{ecuacion}, $R_0( {f}^{2,j}) > R_1, j=N_1+1,\ldots,m.$
\item If $N_1 = N$; that is, the system of equations has full rank, then the system  is solvable if for some specific value of $\nu =0,\ldots,\tau-1$ in \eqref{ecuacion} instead of equating the left hand to zero we equate it to $1$. Doing this for each zero $\zeta$ of $\alpha$ on $\mathcal{C}_1$ and for each $\nu=0,\ldots,\tau-1$, we obtain $N_1$ linearly independent solutions of \eqref{eq:1} which are meromorphic on a neighborhood of $\overline{D_0({\bf f}^1)}$ except for a pole of exact order $\nu+1$ at $\zeta$. On this circle, this would settle the last statement of the theorem. (We will show that on each circle containing more than one zero of $\alpha$ the corresponding system of equations has full rank. This conclusion will be drawn at the very end of the proof of the theorem.)
\end{enumerate}

\medskip

Now let us suppose that $\mathcal{C}_1$ contains only one zero $\zeta_j$ of $\alpha$ of multiplicity $1$; that is $N=1$. If, nevertheless,
\eqref{several}  takes place we could proceed as before, so we will not use this restriction in the arguments that follow. We have,
\[ R_0( {f}^{1,\nu}) \geq R_1 = R_0({\bf f}^1) = |\zeta_j|, \qquad \nu = 1,\ldots,m,\]
with equality for some $\nu$. Without loss of generality we can assume
\[R_0( {f}^{1,1}) = \cdots = R_0( {f}^{1,M}) = R_1, \qquad  1\leq M \leq m,\]
and   $R_0(f^{1,\nu}) > R_1, \nu =M+1,\ldots, m$ (if any).
According to \eqref{cociente} (with $\ell = 1$) and the fact that $\zeta_j$ is the unique zero of $\alpha$ on $\mathcal{C}_1$, we have
\begin{equation}
\label{coc2}
\lim_{n\to \infty} \frac{ {f}^{1,\nu}_{n-1}}{ {f}^{1,\nu}_n} = \zeta_j, \qquad \nu = 1,\ldots,M,
\end{equation}
where $( {f}^{1,\nu}_n)_{n\geq 0}$ denotes the collection of Taylor coefficients of $ {f}^{1,\nu}$, and $\zeta_j$ is a singular point of each $ {f}^{1,\nu}, \nu =1,\ldots,M$. Should $M=1$ we have obtained one solution of \eqref{eq:1}, namely $ {f}^{1,1}$, with radius of convergence $R_1$ and the remaining solutions in ${\bf f}^1$ have radius of convergence larger than $R_1$. We aim to show that if $M > 1$ we can also find one solution of \eqref{eq:1} with radius $R_1$ and  additional $m-1$ linearly independent solutions of \eqref{eq:1} (not necessarily $ {f}^{1,2},\ldots, {f}^{1,m}$) with radius larger than $R_1$.

\medskip

Without loss of generality, we can assume that $ {f}^{1,1}_n \neq 0, n \geq 0$. Indeed, \eqref{coc2} entails that $ {f}^{1,1}_n \neq 0, n \geq n_0$. Let $T_n( {f})$ be the Taylor polynomial of $ {f}$ of degree $n$. Consider the collection of functions $\hat{ {f}}^1,\ldots,\hat{ {f}}^M$, where
\[\hat{ {f}}^\nu(z) = ({ {f}}^{1,\nu}(z) - T_{n_0}({ {f}}^{1,\nu})(z))/z^{n_0}, \qquad \nu =1,\ldots, M.\]
It is easy to verify that these functions are linearly independent, satisfy the recurrence \eqref{eq:1} with the indices $n$ shifted by $n_0$, have radii of convergence equal to $R_1$, have the same singularities as the corresponding $ {f}^{1,\nu}$ on $\mathcal{C}_1$, and $\hat{ {f}}^1_{n} \neq 0, n\geq 0$. Should it be necessary, we derive the desired properties of the functions ${ {f}}^{1,1},\ldots,{ {f}}^{1,M}$ from $\hat{ {f}}^1 ,\ldots,\hat{ {f}}^M$.

\medskip

Let $\lambda = \zeta_j^{-1}$. This point is a root of the characteristic polynomial $p(z)= z^{m}\alpha(1/z)$ of \eqref{eq:1}. Set
\[\gamma_n =  {f}^{1,1}_n/\lambda^n,\qquad n \geq 0. \]
According to Lemma \ref{lem:2}, $(\lambda^n)_{n \geq 0}, \lambda^n =   { {f}^{1,1}_{n}}/{\gamma_n},$ is a solution of \eqref{eq:3}. From \eqref{coc2} we have
\begin{equation}
\label{gamma}
\lim_{n\to \infty} \frac{\gamma_{n+1}}{\gamma_n} = \frac{1}{\lambda}\lim_{n\to \infty} \frac{ {f}^{1,1}_{n+1}}{ {f}^{1,1}_{n}} = 1.
\end{equation}
Consequently, the recurrences \eqref{eq:1} and \eqref{eq:3} have the same characteristic polynomial. Aside from the solution $(\lambda^n)_{n \geq 0}$, \eqref{eq:3} also has the solutions $( {f}^{1,\nu}_{n}/\gamma_n)_{n\geq 0}, \nu = 2,\ldots,M$. Using Lemma \ref{lem:3}, with these solutions of \eqref{eq:3} we can construct $M-1$ linearly independent solutions $(F^{\nu}_{ n})_{n \geq 0}, \nu=2,\ldots,M,$ of \eqref{eq:4} where
\begin{equation}
\label{radii}
F^{\nu}_{n} = \frac{ {f}^{1,\nu}_{n+1}}{\gamma_{n+1}} - \lambda \frac{ {f}^{1,\nu}_{n}}{\gamma_{n}}, \qquad n \geq 0, \qquad \nu=2,\ldots,M.
\end{equation}

\medskip

The radii of convergence of the functions $ {f}^{1,\nu}, \nu =2,\ldots,M$ equal $R_1 = 1/|\lambda|$. This together with  \eqref{gamma} and \eqref{radii} imply that
\[\limsup_{n\to \infty} |F^{\nu}_{n}|^{1/n} \leq 1/R_1, \qquad \nu =2,\ldots,M.\]
According to Buslaev's theorem, for each $\nu = 2, \ldots, M$ the radius of convergence of $F^{\nu}(z) = \sum_{n=0}^{\infty} F^{\nu}_{n} z^n$ is equal to the reciprocal of the absolute value of one of the zeros of the characteristic polynomial $\hat{p}(z)=z^{m-1} + \beta_1 z^{m-1} + \cdots + \beta_{m-1}$ associated with \eqref{eq:4}. According to \eqref{conexion}, and our assumptions, $\lambda$ is not a zero of $\hat{p}$. Therefore,
\[\limsup_{n\to \infty} |F^{\nu}_{n}|^{1/n} := \mu_{\nu} < 1/R_1, \qquad \nu = 2,  \ldots, M. \]

\medskip

Using Lemma \ref{lem:4} and Lemma \ref{lem:3}, with formula \eqref{sol1} we can find $M-1$ linearly independent solutions $\hat{ {f}}^{\nu},\nu=2,\ldots,M$ of the recurrence \eqref{eq:4} whose radius of convergence is greater than $R_1$. Now, using again Lemma \ref{lem:2} with $\gamma_n = \lambda^n/ {f}^1_{1,n}, n\geq 0$ we find a system ${\bf f}^2 = {f^{2,2}, \ldots , f^{2,m}}$ of $m-1$ linearly independent solutions of \eqref{eq:1} each of which has radius of convergence
greater than $R_1 = |\zeta_j|$. Here, $(f^{2,2},\ldots,f^{2,M})$ come from the last application of Lemma \ref{lem:2} whereas $f^{2,\nu} = f^{1,\nu}, \nu = M+1,\ldots,m$. Summarizing, we have found that $R_0(f^{1,1}) = |\zeta_j|$ with $\zeta_j$ a singular point of $f^{1,1}$ and all the functions in ${\bf f}^2$ are linearly independent solutions of \eqref{eq:1} with radius of convergence larger than $|\zeta_j|$.

\medskip

Now we proceed with ${\bf f}^2$ exactly the same way as we did with ${\bf f}^1$ and construct a system ${f}^3$ with $m - N_1 -N_2$ linearly independent solutions of \eqref{eq:1}, where $N_2$ denotes either the rank of the corresponding system of homogeneous linear equations (when the circle $\mathcal{C}_2$ has more than one zero of $\alpha$) or $N_2 = 1$ if $\mathcal{C}_2$ has exactly one zero of $\alpha$. In a finite number of steps, say $r$, we either run out of linearly independent solutions of \eqref{eq:1} because $N_1+\ldots+N_r = m$ or we find at least one non trivial solution of \eqref{eq:1} with radius of convergence $R > |\alpha_m|$.

\medskip

The second possibility is impossible because  according to Buslaev's theorem $R$ has be equal to the absolute value of one of the zeros of $\alpha$. On the other hand, since $N_k, k=1,\ldots,r,$ is less than or equal to the number of zeros, say $m_k$, of $\alpha$ on $\mathcal{C}_k$, if  $N_1+\ldots+N_r = m$, if follows that $N_k = m_k, k=1,\ldots,r$. This happens only when  $\cup_{k=1}^r \mathcal{C}_k$ contains all the zeros of $\alpha$ (we skip no circle at all as we carry out the process). In turn this means that either any given circle contains exactly one zero of $\alpha$ or the rank of the corresponding homogeneous linear system of equations is equal to the number of zeros of $\alpha$ on the circle. As we saw in the proof of the first step this means that associated with each circle $\mathcal{C}_k$ we have $N_k$ linearly independent solutions of \eqref{eq:1} with the properties announced in the statement of the theorem. Since the circles are contained one inside the other the collection of these linearly independent solutions form a fundamental system of solutions of \eqref{eq:1}. \hfill $\Box$

\section{Recurrence relations and Hermite-Pad\'e approximation}

\subsection{Preliminaries.}
Let $\mathbf{f}=\left( f_1,f_2,\ldots,f_d\right)$ be a system of $d$ formal or convergent Taylor expansions about the origin; that is, for each $k=1,\ldots,d$, we have
\begin{equation}\label{serie hp}
f_k(z)=\sum\limits_{n=0}^\infty \phi_{k,n}z^n,\ \ \ \ \phi_{k,n}\in\mathbb{C}.
\end{equation}
In the sequel, ${\bf m} = (m_1,\ldots,m_d) \in \mathbb{Z}^d\setminus {\bf 0}$, where ${\bf 0}\in \mathbb{Z}^d$ denotes the null vector. Set $|{\bf m}| = m_1+\cdots + m_d$.

\begin{Definition}
\label{HP}
Let $({\bf f, m)}$ and $n \geq \max \{m_k: k=1,\ldots,d\}$ be given. Then,  there exist polynomials $q$, $p_k$, $k=1,\ldots,d$, such that
\begin{enumerate}
\item[a.1)] $\deg p_k\leq n-m_k$, $k=1,\ldots,d$, $\deg q\leq |\mathbf{m}|$, $q\not\equiv 0$,
\item[a.2)] $q(z)f_k(z)-p_k=A_kz^{n+1}+\cdots$.
\end{enumerate}
The vector rational function $\mathbf{R}_{n,\mathbf{m}}=\left( p_1/q,\cdots,p_d/q\right)$ is called an $(n,\mathbf{m})$ (type II) Hermite-Pad\'e approximation of $\mathbf{f}$.
\end{Definition}

The existence of such polynomials reduces to solving a homogeneous linear system of $|{\bf m}|$ equations on the $|{\bf m}|+1$ unknown coefficients of $q$. Thus a nontrivial solution exists. Once $q$ is found, the polynomial $p_k, k=1,\ldots,d$ is the Taylor polynomial of degree $n-m_k$ of $qf_k$. Unlike Pad\'e approximants, Hermite-Pad\'e approximants are not uniquely determined (in general). For each $n$ we choose one candidate. Without loss of generality, we take $q$ to have its non zero coefficient of lowest degree equal to $1$. With this normalization, we write $q_{n,{|\bf m|}}$ and $p_{n,{|\bf m|},k}$ instead of $q$ and $p_k$, respectively.

\medskip

Set
\begin{equation}\label{denom}
q_{n,|\mathbf{m}|}(z)=b_{n,|\mathbf{m}|}z^{|\mathbf{m}|}+b_{n,|\mathbf{m}|-1}z^{|\mathbf{m}|-1}+\cdots+b_{n,0}.
\end{equation}
where $b_{n,0} = 1$ if $q(0) \neq 0$ and equals $0$ otherwise.
With this notation, the conditions in Definition \ref{HP} reduce to the following
\[[z^\nu q_{n,|{\bf m}|}f_k]_n = 0, \qquad k=1,\ldots,d,\qquad \nu = 0,\ldots,m_k -1,\]
where $[\,\cdot\,]_n$ denotes the $n$-the Taylor coefficient of $(\cdot)$. In terms of recurrence relations, this means that the system of functions
\begin{equation}
\label{fs}
(f_1,\ldots,z^{m_1-1}f_1,f_2,\ldots,z^{m_2-1}f_2,f_3,\ldots,z^{m_d-1}f_d)
\end{equation}
is made up of solutions of the recurrence relations
\begin{equation}
\label{fse}[  q_{n,|{\bf m}|} f]_n = 0, \qquad n \geq |{\bf m}|,
\end{equation}
of type \eqref{eq:2} (or, what is the same, \eqref{eq:1}) with $\alpha_{n,m}$ replaced by $q_{n,|{\bf m}|}$. We say of type because, we cannot guarantee that $b_{n,|\mathbf{m}|} \neq 0$ or $b_{n,0} = 1$. Aside from this, which will not cause any problem for the application we have in mind (see \eqref{convpol} below), Theorem \ref{teo:1} will allow us to study the analytic properties and detect singularities of the system of functions $\bf f$, under appropriate assumptions on the asymptotic behavior of the sequence of polynomials $q_{n,|{\bf m}|}$.

\medskip

The following notions were introduced in \cite{cacoq2}.

\begin{Definition}
Given $\mathbf{f}=(f_1,\ldots,f_d)$ and $\mathbf{m}=(m_1,\ldots,m_d)\in \mathbb{Z}_+^d\setminus \lbrace \mathbf{0}\rbrace$ we say that $\xi\in\mathbb{C}\setminus\lbrace 0\rbrace$ is a system pole of order $\tau$ of $(\mathbf{f}, \mathbf{m})$ if $\tau$ is the largest positive integer such that for each $s=1,\ldots,\tau$ there exists at least one polynomial combination of the form
\begin{equation}\label{system pole}
\sum\limits_{k=1}^dp_kf_k,\qquad \deg p_k<m_k,\qquad k=1,\ldots,d,
\end{equation}
which is analytic on a neighborhood of $\overline{D}_{|\xi|} = \{z: |z| \leq |\xi|\}$ except for a pole at $z=\xi$ of exact order $s$. If some component $m_k$ equals zero the corresponding polynomial $p_k$ is taken identically equal to zero.
\end{Definition}

\begin{Definition}\label{system sing}
Given $\mathbf{f}=(f_1,\ldots,f_d)$ and $\mathbf{m}=(m_1,\ldots,m_d)\in \mathbb{Z}_+^d\setminus \lbrace \mathbf{0}\rbrace$ we say that $\xi\in\mathbb{C}\setminus\lbrace 0\rbrace$ is a system singularity of $(\mathbf{f}, \mathbf{m})$ if there exists at least one polynomial combination of the form \eqref{system pole} analytic in $D_{|\xi|} = \{z:|z| < |\xi|\}$ and $\xi$ is a singular point of \eqref{system pole}.
\end{Definition}

In this context, the concepts of singular point and pole depend not only on the system of functions $\bf f$ but also on the multi index $\bf m$. For example, poles of the individual functions $f_k$ need not be system poles of $\bf f$ and system poles need not be poles of any of the functions $f_k$ (see interesting examples in \cite{cacoq2}).

\begin{Definition}\label{poly indep}
A vector $\mathbf{f}=(f_1,\ldots,f_d)$ of $d$ formal Taylor expansions is said to be polynomially independent with respect to $\mathbf{m}=(m_1,\ldots,m_d)\in\mathbb{N}^d$ if there do not exist polynomials $p_1,\ldots,p_d$, at least one of which is non-null, such that
\begin{enumerate}
\item[b.1)] $\deg p_k < m_k$, $k=1,\ldots,d$,
\item[b.2)] $\sum\limits_{k=1}^dp_kf_k$ is a polynomial.
\end{enumerate}
\end{Definition}

In particular, polynomial independence implies that for $k=1,\ldots,d$, $f_k$ is not a rational function with at most $m_k-1$ poles and the system of functions \eqref{fs} is linearly independent. Moreover, \eqref{fs} constitutes a fundamental system of solutions of \eqref{fse}.

\subsection{Application of the main result to Hermite-Pad\'{e} approximation.}
In the sequel, we will assume that the sequence \eqref{denom}, $n\geq\max\lbrace m_1,\ldots,m_d\rbrace$ verifies
\begin{equation}
\label{convpol}
\lim\limits_{n\longrightarrow\infty}q_{n,\mathbf{m}}(z)=q_{|\mathbf{m}|}(z) = \prod_{k=1}^{|\mathbf{m}|}(1- z\zeta_k^{-1}),\qquad  \deg q_{|\mathbf{m}|} = |\mathbf{m}|,\qquad q_{|\mathbf{m}|}(0) = 1.
\end{equation}

\begin{Theorem} \label{teo:3}
Suppose that $\mathbf{f}=(f_1,\ldots,f_d)$ is a vector of formal power expansions as in \eqref{serie hp} which is polynomially independent with respect to $\mathbf{m}=(m_1,\ldots,m_d)\in\mathbb{N}^d$. Assume that \eqref{enumeration}, \eqref{limit2} and \eqref{several} take place. Then, each $\zeta_k, k=1,\ldots,|{\bf m}|,$ is a system singularity of $(\mathbf{f}, \mathbf{m})$. Moreover, if $\zeta_k$ verifying \eqref{several} is a zero of multiplicity $\tau_k$ of $q_{|{\bf m}|}$ then it is a system pole of $(\mathbf{f}, \mathbf{m})$ of order $\tau_k$.
\end{Theorem}

\noindent {\bf Proof.}
We have that $\lim_{n\rightarrow\infty}q_{n,\mathbf{m}}=q_{|\mathbf{m}|}$ where $\deg q_{|\mathbf{m}|} = |\mathbf{m}|$ and $q_{|{\bf m}|}(0) = 1$, then there exist $n_0$ such that $b_{n,|\mathbf{m}|}\neq 0, q_{n,0} = 1$, for all $n\geq n_0$.

\medskip

Given $f(z) = \sum_{\nu =0}^{\infty} \phi_\nu z^\nu$ define
\[ \hat{f}(z) := (f(z) - T_{n_0}(f)(z))/z^{n_0 } = \sum_{\nu =0}^{\infty} \phi_{n_0 + \nu} z^{\nu} = \sum_{\nu =0}^{\infty} \hat{\phi}_{\nu} z^{\nu} , \]
where  $T_{n_0}(f)$ is the Taylor polynomial of $f$ of degree $n_0-1$. Notice that the coefficients of $\hat{f}$ are shifted by $n_0$ in relation with the coefficients of $f$; that is
\[\hat{\phi}_\nu = \phi_{n_0 + \nu}, \qquad \nu \geq 0.\]
Consequently,
\[[q_{n_0 + n,|{\bf m}|}f ]_{n_0 +n} = \phi_{n+n_0}+b_{n+n_0,1}\phi_{n+n_0-1}+\cdots+b_{n+n_0,|\mathbf{m}|}\phi_{n+n_0-|\mathbf{m}|}=
\]
\[\hat{\phi}_{n}+b_{n+n_0,1}\hat{\phi}_{n-1}+\cdots+b_{n+n_0,|\mathbf{m}|}\hat{\phi}_{n-|\mathbf{m}|}\]
with $b_{n+n_0,|\mathbf{m}|} \neq 0, n \geq 0$. If we set $\hat{q}_{n,\mathbf{m}}=q_{n+n_0,\mathbf{m}}$, we have
\[[q_{n_0 + n,|{\bf m}|}f ]_{n_0 +n} = 0 \qquad \mbox{if and only if} \qquad [\hat{q}_{n,|{\bf m}|}\hat{f} ]_{n} = 0.\]

Set $\hat{\bf f} = (\hat{f}_1,\ldots,\hat{f}_d)$. Notice that the system of functions
\begin{equation}
\label{sysaux}
(\hat{f}_1,\ldots,z^{m_1-1}\hat{f}_1, f_2,\ldots,z^{m_d-1}\hat{f}_d)
\end{equation}
is linearly independent for otherwise the system $\bf f$ would not be polynomially independent. Consequently, \eqref{sysaux} constitutes a fundamental system of solutions of the recurrence relations
\begin{equation}
\label{ecaux}
[\hat{q}_{n,|{\bf m}|}\hat{f} ]_{n} = 0, \qquad n \geq |{\bf m}|,
\end{equation}
which verifies $\deg(\hat{q}_{n,|{\bf m}|}) = |{\bf m}|, \hat{q}_{n,|{\bf m}|}(0) = 1, n \geq 0,$ and
\[\lim_{n\to \infty} \hat{q}_{n,|{\bf m}|} = q_{|\bf m|}.\]
Therefore, the recurrence relation  \eqref{ecaux} is exactly like \eqref{eq:1}, including the assumption that the coefficient which plays the role of $\alpha_{n,m}$ (which is the coefficient of $q_{n,|{\bf m}|}$ accompanying $z^{|{\bf m}|}$) is different from zero, and fulfills the assumptions of Theorem \ref{teo:2}.

\medskip

Applying Theorem \ref{teo:2} to $(\hat{\bf f},{\bf m})$ we obtain that this pair fulfills   the thesis of Theorem \ref{teo:3}. However, it is easy to verify that $(\bf f,\bf m)$ and $(\hat{\bf f},\bf m)$ have the same system singularities and the same system poles including their order. \hfill $\Box$

\medskip

Two immediate consequences of Theorem \ref{teo:3} which are worth singling out are the following.

\begin{Corollary} \label{cor:1}
Suppose that $\mathbf{f}=(f_1,\ldots,f_d)$ is a vector of formal power expansions as in \eqref{serie hp} which is polynomially independent with respect to $\mathbf{m}=(m_1,\ldots,m_d)\in\mathbb{N}^d$. Assume that all the zeros of $q_{|{\bf m|}}$ verify \eqref{several}. Then, $({\bf f,m})$ has exactly $|{\bf m}|$ system poles which coincide with the zeros of $q_{|{\bf m|}}$, taking account of their order.
\end{Corollary}

Corollary \ref{cor:1} contains the inverse statement, (b) implies (a), in \cite[Theorem 1.4]{cacoq2}.

\begin{Corollary} \label{cor:2}
Suppose that $\mathbf{f}=(f_1,\ldots,f_d)$ is a vector of formal power expansions as in \eqref{serie hp} which is polynomially independent with respect to $\mathbf{m}=(m_1,\ldots,m_d)\in\mathbb{N}^d$. Assume that all the zeros of $q_{|{\bf m|}}$ have distinct absolute value. Then, all the zeros of $q_{|{\bf m}|}$ are system singularities of  $({\bf f,m})$.
\end{Corollary}

This corollary was suggested to take place in the sentence right after \cite[Corollary 5.2]{GY}. In \cite{GY} you may find other results regarding Hermite-Pad\'e approximation which complement this section.

\end{document}